\documentclass{article}

\usepackage{amssymb, graphics}

\newcommand{\R}{\mathbb{R}}
\newcommand{\C}{\mathbb{C}}
\newcommand{\Z}{\mathbb{Z}}
\newcommand{\N}{\mathbb{N}}
\def \T {{\mathcal T}}
\def \P {{P}}
\def \V {{\mathcal V}}
\def \L {{L}}
\def \H {{H}}
\def \a {{\zeta}}
\def \eps {{\epsilon}}

\newtheorem{theorem}{Theorem}

\newtheorem{definition}[theorem]{Definition}

\newtheorem{prediction}[theorem]{Prediction}

\begin{document}
\setcounter{page}{1}

\title
{Critical exponents, conformal invariance
\\
and planar Brownian motion}
\author{Wendelin Werner}

\date {5 July 2000}

\maketitle

\begin{abstract}
In this review paper, we first discuss
 some open problems related to 
two-dimensional self-avoiding paths and critical 
percolation. We then review some closely related 
results (joint work with Greg Lawler and Oded Schramm)
on critical exponents for two-dimensional 
simple random walks, Brownian motions and other
conformally invariant random objects. 
\end{abstract}

\section {Introduction}

The conjecture that the scaling limits of 
many two-dimensional systems in statistical physics 
exhibit conformally invariant behaviour at criticality has
led to striking predictions by theoretical
physicists  concerning, for instance, 
the values of exponents that
describe the behaviour of certian quantities
 near (or at) the critical temperature.
Some of these predictions can be reformulated in 
elementary terms (see for instance the
  conjectures for the number of
self-avoiding walks of length $n$ on a planar lattice).

From a mathematical perspective,
 even if the statement of the conjectures
are clear, the understanding of these predictions and of the 
non-rigorous techniques 
(renormalisation group, conformal field theory,
quantum gravity, the link with
highest-weight representation of some
infinite-dimensional Lie algebras, see 
e.g.  \cite {ID,Cabook}) used by physicists
 has been limited.
Our aim in the present review paper is to present 
some results derived in joint work with Greg Lawler and
Oded Schramm \cite {LW1,LW2,LSW1, LSW2, LSW3,LSWa}
that  proves some of these conjectures, and
improves substantially our understanding of others.
The systems that we will focus on (self-avoiding walks,
critical percolation, simple random walks) correspond in the 
language of conformal field theory
to zero central charge.

We structure this paper as follows:
In order to put our results into
perspective, we start by very briefly describing
 two models (self-avoiding
paths and critical percolation) and some of the 
conjectures that theoretical physicists have produced 
and that are, at present, open mathematical  problems. 
Then, we state theorems derived in joint work with 
Greg Lawler and Oded Schramm \cite {LSW1,LSW2,LSWa}
concerning
critical exponents for simple random walks and planar
Brownian motion (these had been also predicted 
by theoretical physics). We then show how all these 
problems are mathematically related, and,  in particular, why 
the geometry of critical percolation
in its scaling limit,
should be closely related to the geometry of a planar Brownian 
path via a new increasing set-valued process introduced 
by Schramm in \cite {S}.

\section {Review of some prediction of theoretical physics}

\subsection {Predictions for self-avoiding walks}

We first very briefly  describe some  
 predictions of theoretical physics
concerning 
self-avoiding paths in a planar lattice.
For a more detailed mathematical  account on this subject, see 
for instance \cite {MS}.
 
Consider the square lattice $\Z^2$ and define the 
set  $\Omega_n$ of nearest-neighbour paths  of length
$n$ started at the origin that are self-avoiding. 
In other words, $\Omega_n$ is the set of injective 
functions $\{ 0 , \ldots, n \} \to \Z^2$, such 
that $w(0) = (0,0)=0$ and
$| w(1) - w(0) | = \cdots = | w(n) - w(n-1) | = 1 $.
 
The first problem is to understand the asymptotic behaviour of 
the number $a_n := \# \Omega_n$
 of such self-avoiding paths when $n \to \infty$.
A first trivial observation is that for all $n, 
m \ge 1$, $a_{n+m} \le a_n a_m$ because 
the first $n$ steps and the last $m$ steps
 of a $n+m$ long self-avoiding path are self-avoiding paths
of length $n$ and $m$ respectively.
Furthermore,  $a_n \ge 2^n$ because if the path goes
only upward or to the right, then it is self-avoiding. 
This leads immediately to the existence of a constant 
$\mu \in [2,3) $ (called the connectivity
constant of the lattice $\Z^2$)
such that 
$$
\mu:=  \inf_{n \ge 1} (a_n)^{1/n} =
\lim_{n \to \infty} (a_n)^{1/n}
.$$

Note that if one counts the number $a_n'$ 
of self-avoiding paths 
of length $n$ on a triangular lattice, the   
same argument shows that 
$(a_n')^{1/n}$ converges when $n \to \infty$ to 
some limit $\mu' \ge 3$.
The connectivity constant is lattice-dependent.

One can also look at
other regular planar lattices, such as the honeycombe lattice.
We will say that a property is  `lattice-independent' if it 
it holds for all these three `regular' lattices
(square lattice, triangular lattice, honeycombe lattice).
 One possible way to 
describe a larger class
of `regular'  lattices for which  our `lattice-independent'
properties should hold, could   
(but we do not
want to discuss this issue in
detail here) be that they 
are transitive (i.e. for each pair of points, there exists 
a euclidean isometry that maps the lattice onto itself and
one of the two points onto the other) and that rescaled 
simple random walk on this lattice
converges to
planar Brownian motion. For instance, the lattice
$\Z \times 2\Z$ is not allowed.
 
When $L$ denotes such  a 
planar lattice, we denote by $a_{n,L}$ the number 
of self-avoiding paths of length $n$ in the lattice
 starting from 
a fixed point, and by $\mu_{L}
:= \inf_{n \ge 1} (a_{n,L})^{1/n}$ its connectivity 
constant.

A first striking prediction from theoretical 
physics is the following:

\begin {prediction}[Nienhuis \cite {N1}]
For any regular planar lattice ${\L}$,
when $n \to \infty$,
$$
a_n (L)
 =  (\mu_{L})^n n^{11/32 + o(1)}.
$$
\end {prediction}
  
The first important feature is the rational exponent 
$11/32$. The second one is that this result 
does not depend on the lattice 
i.e., the first-order term 
is lattice-dependent while the second is ``universal''.

The following statement is almost equivalent to the 
previous prediction. Suppose that we define under  the
same probability  $P_n$ two independent 
self-avoiding paths $w$ and $w'$ of length  
$n$ on the lattice ${L}$
(the law of $w$ and $w'$ is  the uniform probability on $\Omega_n$).

\begin {prediction}[Intersection exponents version]
When $n \to \infty$,
$$
P_n [ w  \{ 1, 2, \ldots , n\} \cap w'\{ 0,1, \ldots ,n \}
= \emptyset ] 
= n^{-11/32 + o(1)}
. $$
\end {prediction}
  
Indeed,  $w  \{ 1, 2, \ldots , n\}$ and
$w'\{ 0,1, \ldots ,n \}$
are disjoint if and only if the concatenation 
of the  two paths $w$ and $w'$ is a   self-avoiding 
path of length $2n$ so that 
the non-intersection probability is exactly
$a_{2n}/ (a_n)^2$.
 
A second question concerns the typical behaviour of 
a long self-avoiding path, chosen uniformly in $\Omega_n$ when $n$ 
is large.
Let $d(w)$ denote the diameter of 
$w$. 
Theoretical physics 
predicts that the typical diameter is of
order $n^{3/4}$. This had 
already been predicted using a different
(`very non-rigorous') 
argument by Flory \cite {Fl}
in the late 40's.
A formal way to describe this prediction is
the following:

\begin {prediction}[Nienhuis \cite {N1}]
For all $\epsilon >0$
and all regular lattices, 
when $n \to \infty$,
$$
P_n \left[ d (w) \in  [n^{3/4 - \epsilon},
n^{3/4 + \epsilon} ]\right] \to 1 
.$$
\end {prediction}

One of the underlying beliefs that lead to these conjectures 
is that the measure on long  self-avoiding paths,
suitably rescaled, converges when the length goes to infinity,
 towards
a measure on continuous curves,
that posesses some
invariance properties under conformal transformations.
The counterparts of the previous predictions
in terms of this limitting measure then 
go as follows:
Take two independent paths defined under the 
limitting measure, started at distance $\epsilon$
from each other. Then, the probability that the 
two paths are disjoint decays like
$\epsilon^{11/24}$ when $\epsilon \to 0$.
For the second prediction: 
the Hausdorff dimension of a path defined under 
the limitting measure is almost surely $4/3$.

\subsection {Predictions for critical planar percolation}

We now  review some results
predicted by theoretical physics concerning critical planar 
percolation. A more detailed acount on these conjectures 
for mathematicians can be found for instance in \cite {LPS}.
See \cite {G} for a general introduction to percolation. 

Let $p\in (0,1)$ be fixed.
For each edge between
neighbouring points of the lattice, 
erase the edge with probability $1-p$ and 
keep it (and call the edge  open)
with probability $p$ independently 
for all edges. In other words,   
for each edge we toss a biased coin to decide 
whether it is erased or not.
This procedure defines a random subgraph of 
the square grid. It is not difficult to see that 
the large-scale geometry of this subgraph depends
a lot on the value of $p$. 
In particular, there exists a critical value  $p_c$ (called the 
critical probability), such that if $p<p_c$, 
there exists almost surely no unbounded connected component 
in the random subgraph, while if $p>p_c$,
there exists almost surely a unique unbounded connected component
of open edges.
It is 
not very difficult to see 
that the  value $p_c$
of the critical probability is lattice-dependent.
Kesten 
has shown that for $\L=\Z^2$, $p_c = 1/2$.
We are going to be interested in the geometry of large 
connected component when $p = p_c$.

In regular planar  lattices at $p=p_c$ it is  known that 
almost surely  no infinite connected component exists.
However, a simple duality argument shows that in the
square grid, at $p=p_c =1/2$, for 
any  $n \ge 1$, with probability $1/2$, there
exists a path of open edges joining (in that rectangle) 
the bottom and top boundaries of a fixed  $n  \times (n+1)$  rectangle.
This loosely speaking shows that in a big box, with 
large probability, there
exist connected components of diameter 
comparable to the size of the box.
 
Theoretical physics predicts that large-scale
properties of the geometry of critical 
percolation (i.e., percolation on a planar lattice at
its critical probability) are lattice-independent
(even though the value of the critical probability is
lattice-dependent), and, 
in the scaling limit, invariant under conformal transformations;
see e.g. \cite {Ai,LPS}.
Furthermore, physicists have produced  explicit 
formulas that describe some of its features.
 
A first prediction is the following:
Consider critical percolation restricted to an $n \times n$
square (in the square lattice, say), and choose
the connected component $C$ with largest diameter 
(among all connected components). 
The previous observation shows that the diameter of $C$ 
is of the order of magnitude
of $n$.
Define the rescaled discrete outer perimeter of $C$,
 $\partial_n
=\partial C / n$, where $\partial C$ is the boundary of the 
unbounded connected component of the complement of $C$ in the
plane.  

\begin {prediction}[Cluster boundaries \cite {DS,Ca,ADA}]
\label {cb}
The law of $\partial_n$ converges when $n \to \infty$ 
towards a law $\mu$ on continuous paths $\partial$.
Moreover, $\mu$-almost surely: the Hausdorff dimension of  $\partial$ 
is $7/4$, the path $\partial $ is not self-avoiding, and the 
outer boundary of $\partial$ has Hausdorff dimension $4/3$.
\end {prediction}

\noindent
A weaker version of the first part of this prediction is that
for all $\eps>0$,
when $n \to \infty$,
$
\P [ \# \partial C  \in ( n^{7/4 - \epsilon} , n^{7/4 + \epsilon } )  
] \to   1
$.

A second prediction concerns the crossing probabilities 
of a quadrilateral. Suppose that $L>0$ and $l>0$, and
perform critical percolation in the rectangle
$[0, a_n] \times [0, b_n]$ where $a_n$ and $b_n$ 
are the respective integer parts of $Ln$ and $ln$.
Let $x(L,l)$ denote the cross-ratio between the 
four corners of the $L \times l$ rectangle (more precisely,
it is the value $x$ such that there
exists a conformal mapping from the rectangle onto the 
upper half-plane, such that the left and right-hand side 
of the rectangle are mapped onto the intervals
$(-\infty, 0]$ and $[1-x,1]$).

\begin {prediction}[Cardy's formula \cite {Ca2}]
When $n \to \infty$, the probability that there exists a 
path of open edges in the rectangle $[0, a_n]
\times [0,b_n]$ joining  the left  
and right-hand sides of the boundary of the
rectangle
converges  
to  
$$ F(x) =  \frac {3 \Gamma (2/3)}{\Gamma (1/3)^2} x^{1/3}
\null_2F_1 (1/3, 2/3, 4/3; x )$$
where $_2F_1$ is the usual hypergeometric function.
\end {prediction}

\noindent
These results are believed to be lattice-independent.
This second statement  has been predicted by 
Cardy \cite {Ca1,Ca2}, using and generalising 
conformal field 
theory considerations and ideas introduced in 
\cite {BPZ1,BPZ2}. 
Note that this prediction is of a different nature than the
previous ones. It gives an exact formula for
an event in the scaling limit rather than just an exponent.
Assuming conformal invariance, this prediction 
can be reformulated in a half-plane as follows:

\begin {prediction}[Cardy's formula in a half-plane]
For all $a, b >0$, the probability that there exists a  
crossing (a path of open edges) joining
$(- \infty , -an]$ to $[0, bn]$ in the upper half-plane
converges when $n \to \infty$ towards 
$F(b/ (a+b))$. 
 \end {prediction} 

\noindent
Carleson was the first to note
that Cardy's formula takes on a very simple form in an equilateral
triangle.  Suppose $A,B,C$ are the vertices of an equilaterial 
triangle,
say $A=0, B=e^{2i\pi /3}, C = e^{i\pi /3}$.

\begin {prediction}[Cardy's formula in an equilateral triangle]
In the scaling limit (performing 
critical percolation on the grid $\eps \Z^2$, say, and 
letting $\eps \to 0$), 
 the law of the left-most point on  
 $[B,C]$
 that is connected
(in the triangle)
 to  $[A,C]$ is the uniform distribution on $[B,C]$. 
\end {prediction}

\section {Brownian exponents}

We now come to the core of the 
present paper and state some of the results derived 
in the series of papers \cite {LSW1,LSW2,LSW3,LSWa}
concerning critical exponents for planar Brownian
motions and simple random walks.
These are mathematical results (as opposed to the predictions
reviewed in the previous section) that had been  
predicted some 15 or 20 
years ago. The  proofs of these
theorems (that we shall briefly outline in the coming sections)
do 
use conformal invariance, complex analysis and univalent functions.

\subsection {Intersection exponents}

Suppose that $(S_n, n \ge 0)$ and 
$(S_n', n \ge 0)$ are two independent 
simple random walks on the lattice $\Z^2$
that are both started from the origin.
We are interested in the asymptotic 
behaviour (when $n \to \infty$) 
of the  probability that the traces of $S$ and $S'$ 
are disjoint.

\begin {theorem}[\cite {LSW2}]
\label {5/8}
When $n \to \infty$,
$$
\P [ S  \{ 1, 2, \ldots, n \} \cap S' \{0, 1, \ldots, 
n \}  = \emptyset ] 
=
n^{-5/8 + o(1) } 
.$$
\end {theorem}
 
This result had been predicted by Duplantier-Kwon \cite {DK}, see 
also \cite {Dqg,Dcm} for another non-rigorous derivation based
on quantum gravity ideas and  predictions by 
Khnizhnik, Polyakov and Zamolodchikov.
 
Note that, as opposed to self-avoiding walks and
percolation cluster boundaries, the scaling limit of 
planar simple random walk is well-understood mathematically:
It is planar Brownian motion
(this is lattice-independent) and it is invariant 
under conformal transformations (modulo time-change).
This is what makes it possible to prove Theorem \ref {5/8}
as opposed to the analogous prediction for self-avoiding walks.
The scaling limit analog of Theorem \ref {5/8}
is the following:
\begin {theorem}{\cite {LSW2}}
\label {5/8BM} 
Let $B$ and $B'$ denote two independent planar
Brownian motions started at distance 1 from each other.
Then,
when $ t \to \infty$,
$$
\P [ B[0,t] \cap B' [0,t] = \emptyset ] 
= t^{-5/8 + o(1)} 
.$$
\end {theorem}
 
In \cite {LSW1,LSW2,LSW3,LSWa}, analogous results concerning
non-intersection exponents between more than two 
Brownian motions (or simple random walks), in the 
plane or in the half-plane, are derived.

In fact, Theorem \ref {5/8} is a consequence of Theorem \ref {5/8BM}
via an invariance principle argument (see \cite {LP} and
the references therein for the connection between the two 
results).

\subsection {Mandelbrot's conjecture}
 
Let $(B_t, t \ge 0)$ denote a planar Brownian motion.
Define the hull of $B[0,1]$ as the complement 
of the unbounded connected component of $\C \setminus 
B[0,1]$ and define the 
the outer frontier of $B[0,1]$ as the boundary the 
hull of $B[0,1]$.

\begin {theorem}[\cite {LSW1,LSW2,LSWa}]
\label {mandelbrot}
Almost surely, the Hausdorff dimension of the 
outer boundary of $B[0,1]$ is $4/3$.
\end {theorem}

This result had been conjectured by Mandelbrot \cite {M}
based
on simulations and the  analogy with 
the conjectures for self-avoiding walks.
See also \cite {Dcm} for a physics approach based on quantum gravity.
This theorem is in fact a consequence (using  
results  derived by Lawler in \cite {Lhaus})
of the determination of 
the following critical exponent (called disconnection 
exponent):

\begin {theorem}[\cite {LSW1,LSW2,LSWa}]
If $B$ and $B'$ are two independent planar Brownian motions
started from $0$, 
then, when $t \to \infty$,
$$
P [ B[0,t] \cup B' [0,t]
\hbox { does not disconnect } 1 \hbox { from } \infty ]
= t^{-1/3 + o(1)}
.$$
\end {theorem}
Similarly, a consequence of Theorem \ref {5/8BM} is 
that the Hausdorff dimension of the set of cut 
points of the path $B[0,1]$ is almost 
surely $3/4$. Analogously, we get that the
set of pioneer points (i.e., points
$B_t$ that are on the outer boundary of $B[0,t]$)
in a planar Brownian path 
has Hausdorff dimension $7/4$, which -
together with Theorem \ref {mandelbrot} - is reminiscent of 
Prediction \ref {cb}.
 Also, the determination 
of more general exponents give the 
multifractal spectrum of the 
outer frontier of a Brownian path.
 See \cite {Lbuda} and the references therein
(papers by Lawler) for 
the link between critical exponents and 
Hausdorff dimensions. See also \cite {Bef}.

 \section {Universality}

\subsection {In the plane}

A first important step in the proof of  Theorems \ref {5/8BM} and
\ref {mandelbrot} is the  
observation
that some  conformally invariant random probability
measures with a special `restriction' property 
are identical.
Our presentation of universality 
here differs slightly from that  
of the preprints \cite {LW2,LSW1,LSW2,LSW3}. A more 
extended version of the present approach 
and of
its consequences is in preparation.

More precisely, let $D$ denote the open unit disc,
and suppose
that $P^0$ is a  
rotationally invariant probability 
measure defined  on the set of all 
simply connected compact subsets
$K$ of the closed unit disc, that contain the origin and 
such that $K \cap \partial D$ is just one single point
$e(K)$  (we endow this set with a well-chosen $\sigma$-field).
Note that the law of $e(K)$ is the uniform probability
$\lambda$  on 
$\partial D$.

For all $x \in D$, the probability measure $P^x$
is 
defined as the image probability measure of $P^0$
under a M\"obius transformation $\Phi$ 
 from $D$ onto $D$ 
such that $\Phi (0 ) = x$
(rotational invariance of $P^0$ shows that $P^x$
is independent of the actual choice of $\Phi$).
Similarly, for any simply connected open 
set $\Omega \subset \C$ (that is not identical to the whole plane) and
$x \in \Omega$, we define the probability $P^{x,\Omega}$
as the image measure of $P^0$ under a conformal 
transformation $\Phi$ from $D$ onto $\Omega$ with $\Phi ( 0) = x$.
When $\Phi$ does not extend continuously
to the boundary of $D$, one can view
$\Phi (K)$ as the union of $\Phi ( K \setminus \{ e(K)\})$
and the prime end $\Phi ( e(K))$.

\begin {definition}
 We say that $P^0$ is completely conformally invariant 
(in short: CCI) if for any simply connected $\Omega \subset D$, for any
$x \in D$, the measures
$P^{x,D}$ and $P^{x, \Omega}$ are identical when restricted to  the 
family of sets $\{ K  \ : \ K \cap D \subset \Omega \} $.
\end {definition}

Note that (under $P^{x, \Omega}$), $K \cap D \subset \Omega$ 
is equivalent to $\Phi (e(K)) \subset \partial D$.

\begin {theorem}
\label {univ}
There exists a unique conformally invariant probability measure.
It is given by the hull of a Brownian path started from 
$0$ and stopped at its first hitting of the unit circle.
\end {theorem}

\noindent 
{\bf Idea of the proof.}
First, note that the measure $\tilde P$ defined using 
the hull of a stopped  Brownian 
path 
 is indeed CCI, because of 
conformal invariance of planar Brownian motion and the 
strong Markov property.
Then, consider another CCI measure $P$. It is sufficient
to show that for all simply connected $D' \subset D$ that 
contains $0$,
\begin {equation}
\label {above}
P [ K \cap D \subset D'] = \tilde P [ K \cap D \subset D' ] 
\end {equation}
since the family of all such events $\{ K \cap D \subset D'\}$
 is a generating $\pi$-system
of the $\sigma$-field on which we define the 
measure.
Let $\Phi$ denote a conformal map from $D'$ onto 
$D$ with $\Phi(0) = 0 $.
 Then, because of the
CCI property,  both sides of (\ref {above}) are equal to $\lambda (
\Phi ( \partial D \cap \partial D' ))$ and hence equal.

\subsection {In the half-plane}

Analogous (slightly more complicated) 
arguments can be developped 
for subsets $K$ of a domain
that join two parts of the boundary of this domain (as opposed
to joining a point in the
 interior to the boundary, as in the previous 
subsection).
 For convenience, we  consider subsets of 
the equilateral 
triangle $A=0$, $B=e^{2i \pi/3}$, $C= e^{i \pi /3}$.

We  study probability measures $P$ on the 
set of all simply connected 
compact subsets $K$ of $\overline \T$ such that
 $$K \cap \partial \T = [A,A_1] \cup [A, A_2] \cup \{ e(K)\}$$
where $A_1 = A_1(K) \in [A,B]$, $A_2= A_2 (K) \in [A,C]$ 
and $e(K) \in (B,C)$, and such that 
$ \T \setminus K$
consists of exactly two connected components (having
respectively  $[A_1, B] \cup [B, e(K)]$ and
 $[e(K), C] \cup [C,A_2]$ on their boundary).

We say that $P$ satisfies Cardy's formula if the law
of $e(K)$
is the uniform probability measure on $[B,C]$.

We say that $P$ satisfies Property I, if for all 
non-empty $(B', C') \subset [B,C]$, if $\Phi$ denotes
the conformal map from $\T$ onto $\T$ with $\Phi(A)=A$,
$\Phi(B') = B$ and $\Phi (C') = C$, the image measure of $P$
(restricted to $\{e(K) \in (B',C')\})$
under the mapping $\Phi$
is exactly the measure $P$ restricted to $\{A_1 \in [A,\Phi(B)]
\} \cap \{ A_2 \in [A, \Phi (C)]\}$.

Suppose that $K'$ is a  simply connected compact subset  
of $\overline \T$ such that 
$K' \cap [B,C] \not= \emptyset$, $K' \cap ([A,B] \cup [A,C] )
= \emptyset$, and $\T':= \T \setminus K'$ is  simply connected.
Let $\Phi$ denote the conformal mapping from $\T'$ onto $\T$ with 
$\Phi(A)=A$, $\Phi(B)=B$ and $\Phi(C)= C$.
We say that $P$ satisfies Property II, if for all such $K'$,
the image of $P$ (restricted to $\{ K \ : \ K \cap K' = \emptyset \}$)
under $\Phi$ is identical to the 
measure  $P$ restricted to the set $\{ K \ : \ K \cap \Phi
( \partial \T' \setminus \partial \T) = \emptyset \}$.
 
We the say that $P$ is invariant under restriction if it 
satisfies Property I and Property II.
There are other equivalent definitions and formulations of this 
restriction property.

\begin {theorem}
\label {univ2}
There exists a unique probability measure $P$ that satisfies 
Cardy's formula and is invariant under restriction. It is 
given by the hull in $\T$ (see the precise definition 
below) of Brownian motion, started 
from $A$, reflected with oblique angle $\pi/3$ (pointing `away'
from $A$) on [A,B] and $[A,C]$, and stopped
when reaching $[B,C]$.
\end {theorem}

\begin{figure}[h]
\centering
\resizebox{0.5 \textwidth}{!} {\includegraphics {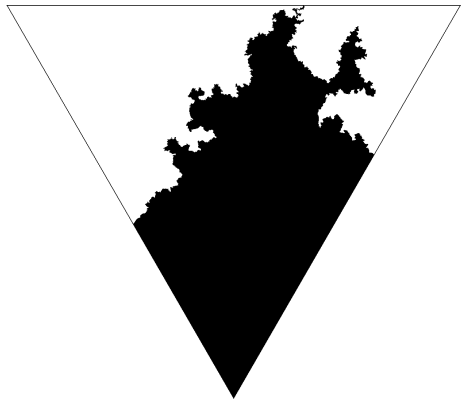}}
\caption{
The $\T$-hull of reflected Brownian motion
}
\end{figure}

Here and in the sequel, the hull
 in the upper half-plane 
(or $H$-hull) of a compact set $K \subset 
\overline H$
with $0 \in K$ is the complement of the 
unbounded connected component of $\overline H \setminus K$.
The hulll in $\T$, or $\T$-hull,  of a compact subset 
of $\overline \T$ that 
contains $A$, is the complement of the union 
of the connected components
of $\overline \T \setminus K$ that have $B$ or $C$ on their 
boundaries.

\medbreak
\noindent
{\bf Idea of the proof.}
Uniqueness follows easily from a similar 
argument as in Theorem \ref {univ}, by identifying 
$P ( K \cap \T \subset  {\T'})$ for a certain class of simply 
connected subsets $\T'$ of $\T$.  
For existence, one needs to verify that 
such a reflected Brownian motion satisfies 
Cardy's formula, and the restriction 
property follows from conformal invariance and the 
strong Markov property for such a reflected Brownian motion.

\medbreak

Let us briefly indicate how this 
reflected Brownian motion is defined (see for instance \cite {VW}
for details).
We first define it in the 
upper half-plane.
Define for any $x \in \R$, 
the vector $u(x) = \exp (i \pi /3)$ if $x \ge0$
and $u(x) = \exp (2 i \pi /3)$ if $x <0$.
Suppose that $B(t)$ is an ordinary planar Brownian path 
started from 
$0$. Then, there exists a unique pair $(Z_t, \ell_t)$ 
of continuous processes such that $Z_t$ takes its values
in $\overline \H$, $\ell_t$ is a non-decreasing real-valued
function with $\ell_0=0$
 that increases only when $Z_t \in \R$, and
$$
Z_t = B_t + \int_0^t u( Z_s) d\ell_s.
$$
The process $(Z_t, t \ge 0)$ is called the reflected 
Brownian motion in $H$ with reflection vector field $u(\cdot)$.
At each time $t$, we define $V_t$ as the $H$-hull
of $Z[0,t]$.
In order to define the $\T$-hull of reflected Brownian motion 
in the triangle that is refered to in the Theorem, consider 
for instance a conformal  mapping $\Phi$ from $\H$
onto the triangle, such that $\Phi(0)= A$,
$\Phi(\infty)= B$, $\Phi (M) =C$ 
for some real $M>0$. Then the $\T$-hull is $\Phi( V_T)$ 
where $T$ is the first time at which $Z$ hits $[M, \infty]$.

\section {Schramm's processes}

\subsection {In the plane}

It is possible to contruct a CCI 
measure $P$ using Loewner's differential equation 
(that encodes a certain class of growing families 
of compact sets) driven by a Brownian motion.
See, for instance, \cite {Dur} for a general introduction to 
Loewner's equation.

For any simply connected compact $K$ such that the
complement of $K$ in the plane is conformally 
equivalent to the complement of a disc,
Riemann's 
mapping theorem shows that there exists a unique $\alpha_K \in \R$
and a unique conformal
map $\hat f_K$
that maps the complement of $K$ onto the complement of 
the unit disk in such a way that 
$$ \hat f_K (z) =  z e^{-\alpha_K}  +  O(1)  $$
when  $z \to \infty$. 

Suppose now that $(\a(t), t \in \R)$ is a continuous function 
taking values on the unit circle.
For any $z \in \C$, define  $f_t (z)$ as the solution 
of the ordinary differential equation 
$$
\partial_t f_t (z) = -
f_t (z) \frac { f_t (z) + \a(t) }{f_t (z) - \a(t)}
$$
such that
$$
\lim_{t \to - \infty}  e^{-t} f_{t} (z) =   z 
$$
For any fixed $z \not= 0$, the mapping $t \mapsto f_t (z)$ is well-defined
up to a possibly infinite  `explosion time'
$T_z$, at which $f_t(z)$ hits
the singularity $\a(t)$ (and we put $T_0 = - \infty$).
Simple considerations show that for any time $t \in \R$, 
$\alpha_{K_t} = t$, and $f_t = \hat f_{K_t}$, where 
$$ K_t = \{ z \in \C \ : \ T_z \le t \}$$
i.e., $f_t$ is the conformal mapping from the 
complement of $K_t$ onto the complement of the unit disc
such that $f_t (z) = z e^{-t} + o(z)$ when $z \to \infty$.

We now take $\a(t) =  \exp ( i \sqrt {\kappa} W_{t} ), t \in \R$,
where $\kappa>0$ is a fixed constant and $(W_t, t \in \R)$
denotes a one-dimensional Brownian motion such that the law
of $\sqrt {\kappa} W(0)$ is the uniform distribution on $[0, 2 \pi]$.
We call $(K_t, t \in \R)$ Schramm's 
radial process with parameter $\kappa$
(in \cite {S,LSW2}, it is refered to as $SLE_\kappa$
for Stochastic Loewner Evolution process).
 Then, let $T$ denote the first time at which $K_t$
intersects the unit circle. 
The set $K_T$ is  a (random) 
simply connected compact set 
that intersects the unit circle at just one point.
 
\begin {theorem}{\cite {LSW1,LSW2}}
\label {SLE6}
If $\kappa=6$, the law of $K_T$ is a CCI probability 
measure.
\end {theorem}

In fact, a more general result (complete conformal invariance of $SLE_6$
as a process) also holds, see \cite {LSW1,LSW2}.
This theorem may seem quite surprising. Note that it 
fails for all other values of $\kappa$.
 The idea of the proof is to 
prove  `invariance'
of the law under infinitesimal deformations of $D$.

A direct consequence of Theorems \ref {SLE6} and \ref {univ} is that 
the law of $K_T$ is identical to the law of the hull of
a stopped Brownian path. This is a rather surprising result
since the two processes (hulls
of planar Brownian motion and $SLE_6$) are 
a priori very different.
For instance, the joint distribution at the first hitting times
 of
the circles of radius 1 and 2 are not the same for $SLE_6$ and 
for the hull of a planar Brownian path.

\subsection {In the half-plane}

Similarly, one can construct natural processes 
of compact subsets of the closed  upper-half plane 
$\overline \H = \{ z \ : \ \Im (z) \ge 0 \}$ that are
`growing
from the boundary'.

For any simply connected compact subset $K$
of $\overline H$ such that $0 \in K$ and $H \setminus K$ is 
simply connected, there exists a unique number
$\beta_K \ge 0$ and  a unique conformal map $\hat g_K$
from $H \setminus K$   
onto $H$ such that $\hat g_K (\infty) = \infty$ and 
$$
\hat g_K (z) = z +  \frac {\beta_K }{ z}  + o \left(
\frac {1}{z} \right)
$$
when $z \to \infty$.
$\beta_K$ which is   increasing in $K$, is 
a half-space analog of capacity.

Suppose now that $(a(t), t \ge 0)$ is a continuous
real-valued function, and define for all $z \in \overline \H$,
the solution $g_t (z)$ to the ordinary differential 
equation
$$
\partial_t g_t (z) = 
\frac { 2 }{ g_t (z) - a(t) }
$$
with $g_0 (z) = z$. This equation is well-defined up
to the (possibly infinite) time 
$
T_z$
at which $g_t(z)$ hits $a(t)$. Then, define
$$
K_t  := \{ z \in \overline \H \ : \ T_z \le t \}.
$$
$(K_t, t \ge 0)$ is an increasing family of 
compact sets, and it is not difficult
to see that for each time  $t\ge 0$, $\beta_{K_t}
= 2t$ and $g_t 
= \hat g_{K_t}$ i.e., $g_t$ is  
the unique conformal mapping from $\H \setminus K_t$
onto $\H$ such that $g_t ( \infty) = \infty$,
and $g_t(z) = z +  2t / z + o(1/z) $ when $z \to \infty$.

When $a(t) = \sqrt {\kappa} W_t$, where $W$ is one-dimensional
Brownian motion with $W(0)=0$, we get a random 
increasing family $(K_t, t \ge 0)$ that we call
Schramm's chordal process with parameter $\kappa$ (it is refered
to as chordal $SLE_\kappa$ 
in \cite {S,LSW1,LSW2,LSW3}).

Note that if we put $Z_t (z) = 
W_t - \kappa^{-1/2} g_t(z)$, then 
$$
Z_t (z) = W_t + \int_0^t \frac {2}{\kappa Z_s(z)} ds 
$$
so that $Z$ can be interpreted as a
(translation of a) complex  
Bessel flow of dimension $1+ (4/ \kappa)$.
If $\kappa<4$, it is easy to see, 
by comparison with a Bessel process,  
that  almost surely $T_z <\infty$ for all $z \in \overline \H$,
in other words, $\cup_{t \ge 0 } K_t = \overline H$.

The scaling properties of
Brownian motion easily show that it is possible to define 
(modulo increasing time-reparametrization) the law
of an increasing process of hulls
in any simply connected set. More precisely, if 
$\Phi$ is a conformal mapping from $\H$ 
onto some simply connected 
domain $\Omega$, 
 then we say that 
$( \Phi (K_t) , t \ge 0)$ is a chordal Schramm process
(with parameter $\kappa$) started from $\Phi(0)$
aiming at $\Phi(\infty)$ in $\Omega$. This process
is well-defined modulo increasing time-reparametrization
(i.e.,  if there exists a (random) continuous increasing $\psi$
such that for all $t$, $K_t= K_{\psi(t)}'$, then we say that the
two processes $K$ and $K'$ are equal).  

For the remainder of this paper, we will
assume that  $\kappa= 6 $ 
as this case  exhibits many very interesting properties.
It is possible to compute explicitely certain probabilities. For
instance, if $a,b>0$,  
$$
P ( T_{-a} < T_b ) 
= F (b/(a+b))
$$
where $F$ is the same as in Cardy's formula 
In fact, much more is true:

\begin {theorem}[\cite {LSW1}]
\label {rest}
Consider $SLE_6$ in the equilateral triangle, started from 
$A$ and aiming at some point $M \in [B,C]$. Let $T$ denote the 
first time $t$ at which $K_t \cap [B,C] \not= \emptyset$.
Then, the law of $K_{T-} = \overline { \cup_{t<T} K_t}$ is 
independent of the choice of $M \in [B,C]$, it satisfies 
Cardy's formula and it is invariant under restriction. 
\end {theorem}

\noindent
Again, this property is only valid when $\kappa = 6$ and there exists a 
more general version in terms of processes \cite {LSW1}.

\subsection {Link with reflected Brownian motion}

From Theorem \ref {rest}, Theorem  \ref {univ2} and the 
strong Markov property,  it follows that 
the chordal process $(K_t, t \ge 0)$ can be reinterpreted in terms
of reflected Brownian motions. In particular, 
the law of $K_{T-}$ is that of the $\T$-hull of 
a stopped Brownian motion with oblique reflection. 
More generally, finite-dimensional marginals of 
the chordal $SLE_6$ process can
for instance  be constructed as follows (other more general
statements with other stopping times hold as well).

Suppose that $J_1, \ldots, J_p$ is a decreasing family  of 
closed subsets of $\overline H$ such that $H \setminus 
J_1, \ldots, H \setminus J_p$ are simply connected.
When $(K_t, t \ge 0)$ is chordal $SLE_6$ in $\overline H$
(started at $0$ and aiming at infinity), we
put, for all $j \le p$,
$$ T_j = \inf \{ t >0 \ : \ K_t \cap J_j \not= \emptyset \}.$$

Let $V$ denote a simply connected compact subset 
of the upper half-plane such that $H'= \H \setminus V$
is simply connected and let $x \in \partial H' $.
Define a reflected Brownian motion $(B_s, s \ge 0)$ 
in $H'$ with oblique reflection 
angle (angle $\pi/3$ on the part of the boundary 
between $x$ and $+\infty$ and reflection $2 \pi / 3$ 
between $-\infty$ and $x$)
as the conformal image in $H'$
of reflected Brownian motion in $H$.
If $J$ is a compact set, define
 the stopping  time 
$S=S (J)$ at which $
V\cup B[0,S]$ intersects $J$ for the first time, and
the hull $\V = \V ( V, x, J)$ of $V \cup B[0, S(J)]$ in $H$.

Now, define recursively (using each time independent 
Brownian motions), $V_0 = \{ 0\}$, $x_0 =  0$,
$$
V_{j+1} = \V ( V_j, x_j, J_j )
\hbox { and }
x_{j+1} = B ( S(J_j )).
$$

\begin {theorem}
\label {approx}
The laws of $(V_{1}, \ldots, V_{p})$
and of $(K_{T_1}, \ldots, K_{T_p})$ are identical.
\end {theorem}

In particular, this implies that the law of $V_p$ is 
identical to that of $K_{T_p}$ and therefore
is independent of 
$p$ and $J_1, \ldots, J_{p-1}$. This can be viewed 
as  a way to
reformulate the restriction property for $SLE_6$ in terms
of reflected Brownian motions only.
So far, there is no direct proof of this fact that does not
use the link with Schramm's process.

\subsection {Relation between radial and chordal $SLE_6$}

Radial and chordal Schramm processes with parameter 6 
are very closely related.
The CCI property and the restriction property in their
`process' versions can 
be  very non-rigorously  described as follows:
The evolution of $K_t$ (for radial and 
chordal $SLE_6$) depends on $K_t$ only
in a local way i.e.,  suppose that at time $t_0$,
$K_t$ is increasing near some point $x \in \partial K_{t_0}$,
then the evolution of $K_t$ immediately after $t_0$ depends only
on how $K_{t_0}$ looks like in the neighbourhood of $x$
(for a rigorous version of this statement, see \cite {LSW1,LSW2}).
The link between radial and
chordal $SLE_6$ (and this is only true when $\kappa=6$)
is that 
 this local evolution is the same for radial and
chordal $SLE_6$, see \cite {LSW2}.
This also leads to a description of the finite-dimensional 
marginal
laws of radial $SLE_6$ in terms
of hulls of  reflected Brownian motion.

\section {Computation of the exponents}

We now give a very brief outline of the proof of 
Theorem \ref {5/8}. 
 As a consequence of the
results stated  in the last 
two sections, hulls of Brownian paths at certain 
stopping
times can be also constructed in an a priori completely
different way using $SLE_6$. 
The latter 
turns out to be much better suited to compute 
probabilities  involving only
the shape of its complement, as $SLE_6$
is  a process that is 
`continuously growing to the outside', whereas  Brownian motion 
 does incursions inside its own hull, so that 
for instance the point at which the Brownian hull `grows' makes
a lot of jumps. 

It is very easy to show that Theorem \ref{5/8}
is a consequence of the following fact that we shall 
now derive:
Let $B$ and $B'$ denote two independent planar Brownian 
motions started from $0$ and $\eps>0$, and
killed when they hit the unit circle. Then, when $\eps \to 0$,
\begin {equation}
\label {5/4}
P [ B \cap B' = \emptyset ] = \eps^{5/4 + o(1) }.
\end {equation}
To derive (\ref {5/4}), note first  that if $h_{K(B)}$
denotes the harmonic measure of $\partial D$ 
at $\eps$ in $D \setminus K(B)$ (here
$K(B)$ is the hull of $B$), then 
$$
P [ B \cap B' = \emptyset | B ] 
= h_{ K(B)}.
$$
Hence, using 
Theorem \ref {univ}, 
shows that 
$$
P [ B \cap B' = \emptyset ]
= 
E [ h_{K(B)} ]
= 
E [ h_{K_T} ]
$$
where $K_T$ is radial $SLE_6$ stopped 
as in Theorem \ref {SLE6}.
The construction
of $K_T$ via Loewner's equation 
makes it possible to study explicitly the 
asymptotic behaviour of $E[ h_{K_T} ]$  
when $\eps \to 0$ by computing 
the highest eigenvalue computation of a differential 
operator, see \cite {LSW2}. From this, 
(\ref {5/4}) follows.

Similar (though more involved) arguments lead  
to Theorem \ref {mandelbrot} and 
to the determination of many other such critical exponents that 
are defined in terms of planar Brownian paths
\cite {LSW1,LSW2,LSW3,LSWa}.

\section {The conjecture for percolation}

Suppose that one performs critical bond percolation ($p=p_c=1/2$) 
in the discrete half-plane $\Z \times \N$. Decide
that all edges of the type $[x,x+1]$ when $x$ is on the real axis 
are erased when $x \ge 0$ and open when $x<0$.
Then, there exists a unique infinite cluster $C_-$ 
formed by the negative half-axis and the union of all 
clusters that are attached to it.
We now explore the outer boundary of this cluster,
starting from the point $(0,0)$. In other words, we 
follow the 
left-most possible  path that 
is not allowed to cross open edges as shown in the picture 
below (for convenience, we draw the line 
that stays at `distance' $1/4$ from $C_-$).
 The open edges are the thin plain lines, and the 
exploration process is the thick plain line.

\begin{figure}[h]
\centering
\resizebox{0.7 \textwidth}{!} {\includegraphics {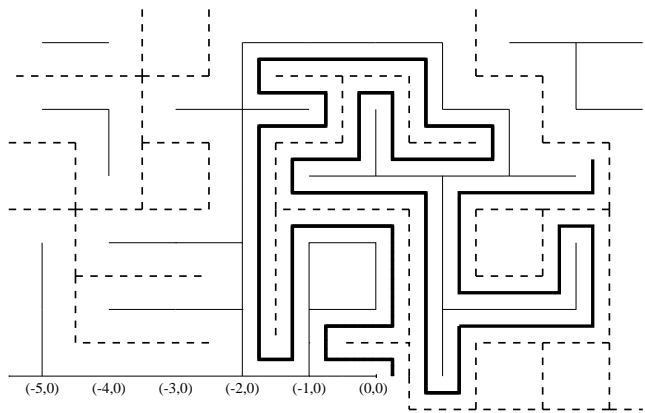}}
\caption{
The discrete exploration process 
}
\end{figure}

Note that this exploration process is almost symmetric
because of the self-duality property of the planar lattice:
It is (almost) the same than exploring the outer boundary of the  
cluster of `open' edges in the dual lattice (i.e.
duals of closed edges in the original 
lattice - the thick dashed edges in the
picture) attached to the `positive half-line'
$(1-i)/ 2 + \N$
(other choices of the lattice make this exploration process
perfectly symmetric).

The conformal invariance 
conjecture  leads naturally
to the conjecture \cite {S,S2} that in the scaling
limit, this  exploration process can be 
described using Loewner's differential equation in 
the upper half-plane, and that 
the driving process $a(t)$ is a continuous symmetric
Markov  process with stationary increments i.e. $a(t)
= \sqrt {\kappa} W_t$ for some parameter $\kappa$. It is 
then easy to identify $\kappa=6$ as the only possible candidate
(looking for instance
at the probability of crossing a square, see \cite {S2}), and the restriction
property
gives additional support to this conjecture.
In particular, note that Cardy's formula for Schramm's
process with parameter $6$ 
would indeed correspond to the crossing probability.

So far, there is no mathematical proof of the 
fact that this exploration process converges to Schramm's
process.
However, Theorem \ref {approx} and the fact that the exploration process 
can be viewed as a discrete random walk reflected on its 
past hull gives  
at least a heuristic hand-waving justification.
 
One can also recover non-rigorously the exponents predicted
for self-avoiding walks (see e.g. \cite {LW2}) using
the exponents derived (rigorously) for Brownian motions
and Schramm processes.

The only discrete models that have been mathematically shown 
to exhibit some conformal invariance properties 
in the scaling limit - apart from simple random walks - 
are those studied by Kenyon in \cite {K1,K2,K3,K4}, namely
loop-erased walks and uniform spanning trees. These are
conjectured to correspond to 
Schramm's processes with parameters 2 and 8, see \cite {S}.
For partial results concerning percolation scaling
limit
and its conformal invariance, see \cite {AB,AB2,BS}.

For some other 
 critical exponents related to Hausdorff dimensions
of conformally invariant exceptional subsets of the planar 
Brownian curve (such as pivoting cut points for
instance) that we do not know (yet?) the value of, see   
\cite {Bef}.

\medbreak
\noindent
{\bf Acknowledgements.}
Without my coauthors 
Greg Lawler and Oded Schramm, 
this paper would of course not exist.
I also use this opportunity to 
thank my Orsay  colleagues Richard Kenyon and
Yves Le Jan for many very stimulating and inspiring discussions, as
well as Vincent Beffara for the nice picture (Figure 1).

\itemsep=\smallskipamount

Laboratoire de Math\'ematiques,

Universit\'e Paris-Sud,

B\^at. 425

91405 Orsay cedex, France

wendelin.werner@math.u-psud.fr

\end{document}